\documentclass[letterpaper, 10 pt, conference]{ieeeconf}  

\IEEEoverridecommandlockouts                              

\usepackage{amsmath} 
\usepackage{amssymb}  
\usepackage{algorithmic}
\usepackage{mathtools}

\newcommand*\diff{\mathop{}\!\mathrm{d}}
\newcommand\norm[1]{\left\lVert#1\right\rVert}
\usepackage[ruled,vlined]{algorithm2e}
\usepackage{comment}
\usepackage{color,soul}
\usepackage{mwe}
\usepackage{svg}
\usepackage{hyperref}
\usepackage{balance}
\usepackage{lipsum} 

\usepackage{psfrag}
\usepackage{pstool}
\graphicspath{{Figures/}}

\title{\LARGE \bf
Solving optimal control problems with non-smooth solutions using an integrated residual method and flexible mesh
}

\author{Lucian Nita$^{1}$, Eric C. Kerrigan$^{1}$, Eduardo M. G. Vila$^{1}$ and Yuanbo Nie$^{2}$
\thanks{This work has received funding from the Engineering and Physical Sciences Research Council under a Doctoral Training Grant (reference number: EP/T51780X/1)}
\thanks{$^{1}$ Department of Electrical \& Electronic Engineering, Imperial College London, SW7~2AZ London, UK, email: 
        {\tt\small \{lucian.nita16, e.kerrigan, emg216\}@imperial.ac.uk}}%
\thanks{$^{2}$ Department of Automatic Control and Systems Engineering, University of Sheffield, S1~3JD Sheffield, UK, email:
        {\tt\small y.nie@sheffield.ac.uk}}%
}

\begin{document}

\maketitle
\thispagestyle{empty}
\pagestyle{empty}

\begin{abstract}
Solutions to optimal control problems can be discontinuous, even if all the functionals defining the problem are smooth. This can cause difficulties when numerically computing  solutions to these problems. While conventional numerical methods assume state and input trajectories are continuous and differentiable or smooth, our method is able to capture discontinuities in the solution by introducing time-mesh nodes as decision variables. This allows one to obtain a higher accuracy solution for the same number of mesh nodes compared to a fixed time-mesh approach. Furthermore, we propose to first solve a sequence of suitably-defined least-squares  problems to ensure that the error in the dynamic equation is below a given tolerance. The cost functional is then minimized subject to an inequality constraint on the dynamic equation residual. We demonstrate our implementation on an optimal control problem that has a chattering solution. Solving such a problem is difficult, since the solution involves infinitely many switches of decreasing duration. This simulation shows how the flexible mesh is able to capture discontinuities present in the solution and achieve superlinear convergence as the number of mesh intervals is increased.
\end{abstract}


\section{INTRODUCTION}
\label{sec:section1}
Solving a sequence of constrained optimal control problems (OCPs) in real-time is a very powerful technique typically used in model predictive control (MPC). However, it can also be challenging in practice to reliably compute a solution since the continuous-time, infinite-dimensional OCP is solved using finite-dimensional numerical solvers. Since we are solving a discretized version of the original problem, the obtained solution is not guaranteed to be feasible for the original continuous-time OCP. Additionally, it can also be difficult to preserve the accuracy of the solution and  obtain superlinear convergence as the number of mesh nodes is increased, especially in problems with discontinuous solutions. 

The most commonly used direct transcription method for solving optimal control problems is direct collocation, which is considered to be the current state of the art~\cite{rawlings2017model}. While collocation has the advantage of being able to handle complex dynamical models, collocation has the fundamental drawback of not guaranteeing an acceptable accuracy in between the collocation points. 

The idea of using integrated residuals as part of the transcription process overcomes some of the limitations of collocation~\cite{nie2019efficient,nie2022solving}. Compared to classical time-marching schemes (shooting methods) or point-wise residual minimisation (collocation), integrated residual methods have the benefit of producing a solution with a more uniform error over the whole time domain. 

State-of-the-art mesh refinement methods are hp-adaptive methods~\cite{liu2015adaptive, paiva2015adaptive}. Some advanced methods have discontinuity detection schemes, but most current refinement strategies rely on knowing beforehand whether the solution will be discontinuous and cannot provide an efficiency comparable to the continuous case for general problems. By adding a flexible mesh, as in this paper, it is possible to develop a method that has similar convergence properties for problems with discontinuous or continuous solutions.


This paper extends the work in~\cite{ECCFlexibleMesh}, which used the integrated residual method for solving differential equations, feasibility problems and constraint satisfaction problems with a flexible mesh. The central contribution of this paper is to extend such integrated residual methods to the solution of OCPs. In the transcription process, a flexible time-mesh will be introduced in order to achieve superlinear convergence for discontinuous problems during the mesh refinement phase. Moreover, the proposed algorithm is able to solve difficult problems to a user-defined accuracy. A numerical example shows how our method performs on a control problem with an optimal chattering solution.

In Section~\ref{sec:section2}, we introduce the optimal control problem formulation. In Section~\ref{sec:section3A} we present the integrated residual method for transcribing the constrained optimal control problem.  Section~\ref{sec:section3B} describes the concept of a flexible time mesh and how this method can improve convergence. Section~\ref{sec:section3C} presents an algorithm for solving an OCP to a user-specified accuracy. Section~\ref{sec:section4} demonstrates how the proposed algorithm from Section~\ref{sec:section3C} can be used jointly with a flexible mesh scheme to solve an optimal control problem and construct a Pareto front between solution accuracy and a lower bound on the optimal cost. Section~\ref{sec:section5} provides a summary of the main findings presented in this paper and discusses potential improvements and future works.

\section{Problem definition}
\label{sec:section2}
The objective functional of many optimal control and estimation problems can be written in the general Bolza form
\begin{subequations}
    \label{eq:equation1}
    \begin{align}
        \min_{x(\cdot),u(\cdot)} & \phi(x(t_0),x(t_f),t_0,t_f) + \int_{t_0}^{t_f}{L(x(t),u(t),t)} \diff t \label{eq:equation1A}\\
        \textrm{s.t.} \quad & F(\dot{x}(t),x(t),u(t),t)=0 \quad \forall t \in[t_0,t_f], \label{eq:equation1B}\\
        & G(\dot{x}(t), x(t), u(t), t) \leq 0 \quad \forall t \in [t_0,t_f], \label{eq:equation1C}
    \end{align}
where  $x : \mathbb{R} \rightarrow \mathbb{R}^{N_x}$ are the state variables and constrained to be continuous, $\dot{x} : \mathbb{R} \rightarrow \mathbb{R}^{N_x}$  are the time derivatives of the state $x$, and $u : \mathbb{R} \rightarrow \mathbb{R}^{N_u}$ are the control inputs. The function $F : \mathbb{R}^{N_x} \times \mathbb{R}^{N_x} \times \mathbb{R}^{N_u} \times \mathbb{R} \rightarrow \mathbb{R}^{N_F}$, which contains the dynamical model of the system, defines a set of $N_F$ equality constraints that have to be satisfied by the controlled system. $G : \mathbb{R}^{N_x} \times \mathbb{R}^{N_x} \times \mathbb{R}^{N_u} \times \mathbb{R} \rightarrow \mathbb{R}^{N_g}$ defines $N_g$ path inequality constraints. $\phi:\mathbb{R}^{N_x} \times \mathbb{R}^{N_x} \times \mathbb{R} \times \mathbb{R} \rightarrow \mathbb{R}$ is the Mayer cost functional, also called the boundary cost, with $t_0\in\mathbb{R}$ and $t_f\in\mathbb{R}$ being the initial and final times, respectively. $L:\mathbb{R}^{N_x} \times \mathbb{R}^{N_u} \times \mathbb{R} \rightarrow \mathbb{R}$ is the Lagrange cost functional, also known as the path cost. Additionally, the problem may have one or more boundary constraints of the form
\begin{align}
    \Psi_E(x(t_0), x(t_f), t_0, t_f) & = 0, \label{eq:equation1D}\\
    \Psi_I(x(t_0), x(t_f), t_0, t_f) & \leq 0, 
    \label{eq:equation1E}
\end{align}
\end{subequations}
where $\Psi_E : \mathbb{R}^{N_x} \times \mathbb{R}^{N_x} \times \mathbb{R} \times \mathbb{R} \rightarrow \mathbb{R}^{N_E}$ are the boundary equality constraints, and $\Psi_I : \mathbb{R}^{N_x} \times \mathbb{R}^{N_x} \times \mathbb{R}  \times \mathbb{R} \rightarrow \mathbb{R}^{N_I}$ are the boundary inequality constraints.

\section{Solution method}
\label{sec:section3}
In most real-time control applications we are heavily constrained by the computational time. As a result, solving the entire problem~\eqref{eq:equation1} using direct collocation has two fundamental drawbacks. Firstly, the designer is not able to control the solution accuracy over the entire time interval without a posteriori computing the error and conducting mesh refinement procedures. Consequently, existing state-of-the-art methods may fail to ensure constraint satisfaction. Secondly, existing schemes cannot terminate early and return the best feasible solution that was achieved in the given amount of computational time. 
Since one often wants to focus on fast constraint satisfaction, we propose to initially solve a feasibility problem and refine the mesh until the  dynamic constraints are satisfied to a given accuracy. We will then use the obtained solution as an initial guess to the optimal control solver, which optimizes a transcribed version of the original problem \eqref{eq:equation1}. 

\subsection{Integrated residual transcription}
\label{sec:section3A}
In the transcription process the infinite-dimensional OCP~\eqref{eq:equation1} has to be converted into a finite-dimensional nonlinear programming problem (NLP). In order to achieve this, the state $x(\cdot)$ and input $u(\cdot)$ trajectories need to be parametrized by a finite number of decision variables $s_{i}^{j}$ and $c_{i}^{j}$ where the subscript $i$ denotes the interval number and $j$ denotes the index of the nodal point in interval $i$, as will be described later. Using a linear combination of these decision variables, approximating functions $\tilde{x}: \mathbb{R} \rightarrow \mathbb{R}^{N_x}$ and $\tilde{u}: \mathbb{R} \rightarrow \mathbb{R}^{N_u}$ can be constructed.     

Before aiming to minimize the objective, in most applications it is critical to ensure the constraints are satisfied to a user-defined accuracy. For this purpose we will introduce an error metric $\epsilon_R\in\mathbb{R}$ defined as
\begin{equation}
    \epsilon_R:=\frac{1}{(t_f-t_0)\cdot N_F}\int_{t_0}^{t_f}{\norm{ F(\dot{\tilde{x}}(t),\tilde{x}(t),\tilde{u}(t),t)}^2_2}\diff t
    \label{eq:equation2}
\end{equation}
based on the integral of the 2-norm squared of the dynamic equation residual. The residual $\norm{ F(\dot{\tilde{x}}(t),\tilde{x}(t),\tilde{u}(t),t)}^2_2$ indicates how well the numerical solution $(\tilde{x}(\cdot),\tilde{u}(\cdot))$ satisfies the dynamic constraint~\eqref{eq:equation1B} over the whole time interval. 
In contrast with direct collocation that enforces constraint~\eqref{eq:equation1B} exactly, but only at a finite number of nodes called collocation points, our method uses quadrature rules to integrate the residual over the whole interval $[t_0, t_f]$, thus guaranteeing a certain level of accuracy in between the collocation points. Note also that the above error metric is a scaled version of the integrated residual where the scaling factor $\frac{1}{(t_f-t_0)\cdot N_F}$ is introduced to average out the residual over the interval $[t_0,t_f]$ and over all components of the dynamics function~$F$. 

Lagrange polynomial basis functions are often used to express the approximating functions $\tilde{x}$ and $\tilde{u}$~\cite[Sect.~1.17.1]{betts2010practical}. The possible solution space is defined by the basis functions used to represent approximation functions $(\tilde{x}(\cdot),\tilde{u}(\cdot))$. As a consequence, the exact solution $(x(\cdot),u(\cdot))$ may not be representable in that solution space, which implies that an exact representation of the constraint~\eqref{eq:equation1B} can never be achieved in finite time (the integrated residual $\epsilon_R$ can asymptotically converge to zero only in the limit as the time-mesh is refined and the number of discretization points is increased).

To reduce the approximation error (as quantified by $\epsilon_R$) there are two fundamental refinement strategies: \begin{itemize}
    \item h-refinement involves splitting the entire time domain  $[t_0, t_f]$ into $N$ subdomains, i.e.\ subintervals $[t_i,t_{i+1}]$ such that
    \begin{subequations}
    \label{eq:equation3}
    \begin{align}
    &\mathcal{T}_i:= [t_{i}, t_{i+1}]  \subset[t_0,t_f], \ \forall i \in \{0, \dots,N-1\},\label{eq:equation3A}\\
    &\cup_{i=0}^{N-1} [t_i,t_{i+1}]=[t_0,t_f],\label{eq:equation3B}\\
    &
    t_i < t_{i+1},\ \forall i \in \{0, \dots,N-1\}
    \label{eq:equation3C} 
    \end{align}
    \end{subequations}
    where $t_{N}=t_{f}$. The refinement variable is therefore the number of subdomains $N$. 
    \item p-refinement relies on constructing a polynomial approximation $\chi_i$ of degree $a$ inside each subdomain $[t_i,t_{i+1}]$ such that for all $i\in\{0, \dots,N-1\}$:
    \begin{subequations}
    \label{eq:equation4}
    \begin{align}
        \tilde{x}(t)&:=\chi_i(t),\ \forall t\in[t_i,t_{i+1}]
        \label{eq:equation4A}\\
        \chi_i(t)&=\frac{\sum_{j=0}^{a}{\frac{w_i^j}{t-\tau_i^j}\cdot s_i^j}}{\sum_{j=0}^{a}\frac{w_i^j}{t-\tau_i^j}\label{eq:equation4B}},
    \end{align}
    \end{subequations}
    where $s_i^j=\chi_i(\tau_i^j)=\tilde{x}(\tau_i^j)$ are NLP decision variables, $w_i^j$ are polynomial weights and $\tau_i^j$ are polynomial nodes~\cite{berrut2004barycentric}. In this case, polynomial refinement means increasing the polynomial degree $a$. Note a similar expression for $\tilde{u}(\cdot)$ can be derived with $b$ denoting the polynomial degree of $\tilde{u}(\cdot)$.  
\end{itemize}
Note that these elementary methods can both be used during the mesh refinement process leading to the so-called hp-type refinement method. 

To enforce state continuity at mesh nodes $t_i$, the additional constraints
\begin{align}
    \tilde{\chi}_i(t_{i+1})&=\tilde{\chi}_{i+1}(t_{i+1}),\quad \forall i\in \{0,\dots,N-2 \},\label{eq:equation5}
\end{align}
are enforced by using the same variable $s_i^a=s_{i+1}^0$ to represent both $\tilde{\chi}_i(t_{i+1})$ and $\tilde{\chi}_{i+1}(t_{i+1}), \forall i \in \{0,\dots,N-2\}$.

\subsection{Residual minimization problem: Improving accuracy to ensure feasibility}
\label{sec:section3B}
To efficiently solve feasibility and control problems with discontinuous solutions, which are otherwise difficult to solve, we will use an integrated residual method to tackle the dynamic constraints. The idea is similar to what \cite{eason1976review}, \cite{mortari2019high} have proposed for solving differential equations and what has been used in \cite{ECCFlexibleMesh} for solving dynamic feasibility problems. 

The first step of our approach is to solve a feasibility problem that aims to satisfy constraints~\eqref{eq:equation1B}--\eqref{eq:equation1E} to a given tolerance. This feasibility problem is converted into a minimisation problem that minimizes the integrated residual of the dynamics model $\epsilon_R$ below a user-specified value. 

In numerical simulations, integrals from~\eqref{eq:equation1A} and~\eqref{eq:equation2} have to be approximated using a $Q$-point Gaussian quadrature rule. Since $F$ is a general nonlinear function, the approximation of the above integrals will not be exact. Apart from the residual error $\epsilon_R$ appearing as a result of the discretization, another numerical error is introduced, namely the quadrature error
\begin{equation}
    \epsilon_Q\coloneqq\Big|\epsilon_R-\sum_{i=0}^{N-1}{\sum_{k=1}^{Q}{\sigma_i^k \cdot \norm{ F(\dot{\tilde{x}}(\rho_i^{k}),\tilde{x}(\rho_i^{k}),\tilde{u}(\rho_i^{k}),\rho_i^{k})}^2_2\Big|}}
    \label{eq:equation6}
\end{equation} where $\sigma_i^k$ for $k\in\{1,\dots,Q\}$ are the $Q$ quadrature weights associated with the integration interval $[t_i,t_{i+1}]$, appropriately scaled by~$N_F$ and interval length to include the initial factors in~\eqref{eq:equation2}, while~$\rho_i^{k}$ are the quadrature nodes for the interval $[t_i,t_{i+1}]$. 

In order to validate the obtained solution, we need to check whether $\epsilon_Q$ is sufficiently small by recomputing the integrals with a higher quadrature order. If the difference between the new value and the solution obtained from the optimization problem is above a certain tolerance $\varepsilon_{quad,tol}$, the problem needs to be resolved using a higher value for $Q$. 

We rely on mesh refinement to select appropriate values for $N$, $a$ and $b$. Note however that conventional mesh refinement strategies applied to a fixed time-mesh with nodes at predefined locations may not always achieve superlinear converge to the solution as the number of nodes is increased. Consider for example the case when a discontinuity in the solution $u(\cdot)$ is located in the interval $(t_{i},t_{i+1})$. In this case, a numerical approximation of this discontinuous function is obtained using a continuous polynomial basis (as described in Section~\ref{sec:section3A}). In general, unless a mesh node is located exactly at the point of discontinuity, a Gibbs phenomenon can occur when interpolating a discontinuous function with a continuous one. This leads to interpolation overshoots that cannot  be eliminated in general by mesh refinement schemes and will cause the error to plateau and not decrease beyond a certain level. 

In order to achieve superlinear convergence in cases where discontinuities are present and produce an accurate solution, we propose including mesh points $t_i$ as decision variables in the NLP formulation. As a result,  time nodes are allowed to move towards regions non-smoothness.

Recall that standard direct collocation methods compute an integral of the residual  only after the NLP has been solved~\cite{betts2010practical}; they do not directly constrain the integral of the residual while computing a solution to the NLP. 
It follows that introducing mesh nodes as decision variables in standard collocation methods  can result in less accurate solutions than those with fixed nodes, unless care is taken. This argument also motivates the interdependence between the use of a flexible mesh and the integrated residual transcription method proposed here.

In an ideal scenario where no quadrature error is present, nodes can be allowed to move freely in the domain according to \eqref{eq:equation3}. However, since quadrature error increases as the intervals expand, we still need to constrain the allowed flexibility of the nodes. For a fixed parameter $\phi \in [0,1)$ we impose upper and lower bounds on the interval length
\begin{subequations}
 \label{eq:equation7}
    \begin{align}
        t_{i+1}-t_i & \leq (1+\phi)\cdot\frac{t_f-t_0}{N},\   \forall i\in\{1,\dots,N-1\},
    \label{eq:equation7A}\\
    t_{i+1}-t_i & \geq (1-\phi)\cdot\frac{t_f-t_0}{N},\   \forall i\in\{1,\dots,N-1\}.
            \label{eq:equation7B}
    \end{align}
\end{subequations}
Since $\phi \in [0,1)$ by definition, the order of nodes is preserved and intervals do not overlap, as required by~\eqref{eq:equation3C}.

For a better sparsity structure of the Hessian along with the ease of implementing continuity constraints \eqref{eq:equation5} and boundary constraints \eqref{eq:equation1D}, \eqref{eq:equation1E} the decision vector $z \in \mathbb{R}^{N\cdot(N_x\cdot(a)+N_u\cdot(b+1)+1)+N_x+1}$ is ordered as 
\begin{equation}
    z:=(s_0^0, t_0, s_0^1,\dots,s_0^{a-1},c_0^0,\dots,c_0^b,s_0^a,t_1,s_1^1\dots,s_{N-1}^a,t_N).
    \label{eq:equation8}
\end{equation}

Hence, the minimum residual solution along with the optimal node locations can be computed from the optimization problem
\begin{subequations}
\label{eq:equation9}
    \begin{align}
    \min_{z} \quad & \sum_{i=0}^{N-1}{\sum_{k=1}^{Q}{\sigma_i^k \cdot \norm{ F(\dot{\tilde{x}}(\rho_i^{k}),\tilde{x}(\rho_i^{k}),\tilde{u}(\rho_i^{k}),\rho_i^k)}^2_2}}\label{eq:equation9A}\\
    \textrm{s.t.} \quad & 
    G(\dot{\tilde{x}}(\tau_i^j), \tilde{x}(\tau_i^j), \tilde{u}(\tau_i^j), \tau_i^j) \leq 0 \quad \forall \tau_i^j\in\tau \label{eq:equation9B},\\
    & \Psi_E(s_0^0, s_{N-1}^a, t_0, t_f) = 0, \label{eq:equation9C}\\
    &\Psi_I(s_0^0, s_{N-1}^a, t_0, t_f) \leq 0,\label{eq:equation9D}\\
    &(1-\phi)\frac{t_f-t_0}{N} \leq t_{i+1}-t_i \leq (1+\phi)\frac{t_f-t_0}{N}, \label{eq:equation9E}\\
    &\tilde{\chi}_i(t_{i+1})=\tilde{\chi}_{i+1}(t_{i+1})\quad \forall i\in \{0,\dots,N-2 \},\label{eq:equation9F}
    \end{align}
\end{subequations}
where path inequality constraints \eqref{eq:equation9B} are implemented at the support time points $\tau_i^j$.
Note that, since time-mesh nodes are added in the decision vector, quadrature points $\rho_i^k$ and internal supports $\tau_i^j$ become functions of $t_i$ and $t_{i+1}$. Hence, these values need to be shifted and scaled accordingly. 
\subsection{Cost minimization problem: From constrained control to optimal control}
\label{sec:section3C}
In this paper we will  focus on h-refinement. 
Starting from a coarse mesh (small $N$) problem~\eqref{eq:equation9} is solved repeatedly until a desirable user-defined tolerance $\epsilon_{tol}$ on~\eqref{eq:equation2} has been reached. Even if the convergence is superlinear with respect to the number of subdomains $N$, the performance can further be improved by providing a good initial guess obtained from interpolating the solution obtained from the previous solution with a coarser mesh. 

After the desired tolerance has been reached,  the mesh parameters $N$ and $a$ are fixed and  the cost functional is minimized by solving 
\begin{subequations}
    \label{eq:equation11}
    \begin{align}
        \min_{z} \, & \phi(x(t_0),x(t_f),t_0,t_f) + \sum_{i=0}^{N-1}{\sum_{k=1}^{Q}{\sigma_i^k L(x(\rho_i^k),u(\rho_i^k),\rho_i^k)}} \label{eq:equation11A}\\
        \textrm{s.t.} \, &  \sum_{i=0}^{N-1}{\sum_{k=1}^{Q}{\sigma_i^k \cdot \norm{ F(\dot{\tilde{x}}(\rho_i^{k}),\tilde{x}(\rho_i^{k}),\tilde{u}(\rho_i^{k}),\rho_i^k)}^2_2}} \leq \epsilon_{tol}\label{eq:equation11B}\\
        &\eqref{eq:equation9B}, \eqref{eq:equation9C}, \eqref{eq:equation9D}, \eqref{eq:equation9E}, \eqref{eq:equation9F},
        \label{eq:equation11C}
    \end{align}
\end{subequations}

Since we have successfully solved \eqref{eq:equation9}, we know the dynamics constraint \eqref{eq:equation11B} should be feasible. We also have an upper bound on the cost and a sufficiently good initial guess needed to efficiently warm-start problem \eqref{eq:equation11}. 

The proposed algorithm for solving optimal control problems using integrated residual transcription method is outlined in Algorithm \ref{alg:algorithm1}. 
\begin{algorithm}[tb]
 \caption{Algorithm for solving OCPs in form~\eqref{eq:equation1}}
 \label{alg:algorithm1}
 \begin{algorithmic}[1]
\REQUIRE $\epsilon_{tol}$, $\varepsilon_{quad,tol}$ and initial values for $N$, $a$, $b$, $Q$ and $z^*$.
\REPEAT
    \STATE $z_0 \gets z^*$
    \STATE $z^* \gets \arg\min_{z} \eqref{eq:equation9A}$ s.t. \eqref{eq:equation9B}, \eqref{eq:equation9C}, \eqref{eq:equation9D}, \eqref{eq:equation9E}, \eqref{eq:equation9F}
    \IF{$\epsilon_Q \leq \varepsilon_{quad,tol}$}
        \STATE $\epsilon_R \gets \min_{z} \eqref{eq:equation9A}$ s.t. \eqref{eq:equation9B}, \eqref{eq:equation9C}, \eqref{eq:equation9D}, \eqref{eq:equation9E}, \eqref{eq:equation9F}
        \STATE increase $N$
    \ELSE
        \STATE increase $Q$
    \ENDIF
\UNTIL{$\epsilon_R \leq \epsilon_{tol}$}

\STATE $z_0 \gets z^*$
\STATE $z^* \gets \min_{z} \eqref{eq:equation11A}$ s.t. \eqref{eq:equation11B}, \eqref{eq:equation11C}
\STATE $\tilde{x}, \tilde{u} \gets \operatorname{interpolate}(z^*)$
 \end{algorithmic} 
 \end{algorithm}
In the initialization phase, mesh variables $N$, $a$ and $b$ are set to small values. The integrated residual minimization problem \eqref{eq:equation9} is then solved using warm starting and the solution checked if the obtained residual is below the threshold tolerance $\epsilon_{tol}$. Note on line~2 of Algorithm~\ref{alg:algorithm1} the pseudo-code notation is simplistic, but the initial guess $z_0$ is not directly set to $z^*$ since the size of $z$ increases as the mesh is refined. However, $z_0$ will be an expanded and interpolated version of $z^*$ as previously explained. Any suitable method can be used for increasing $Q$ and $N$; we chose to double  $N$ and $Q$ every time in our example in Section~\ref{sec:section4}. Once problem \eqref{eq:equation9} has been solved and the residual  minimized,  problem \eqref{eq:equation11}, which is a transcribed version of problem \eqref{eq:equation1}, is solved using available NLP solvers.

\section{Numerical Results}
\label{sec:section4}
An optimal control problem solver based on the integrated residual method was developed in the Julia~v1.6 programming language. The package makes use of barycentric interpolation routines as described in~\cite{berrut2004barycentric} to parametrize the state and input variables. Numerical integration was performed using Gaussian quadrature as detailed in~\cite{quadgk}. Derivative information was obtained using automatic differentiation (AD) tools~\cite{innes2018don} and supplied to the solver as the gradient and Hessian of the Lagrangian function.  
The solver includes an implementation of the flexible mesh scheme of Section~\ref{sec:section3A} along with a fixed mesh version, where time nodes $t_i$ are chosen to be in predefined locations and  not included  as decision variables. In our implementation, Chebyshev type~2 interpolation nodes and weights were used, since we want internal supports $\tau_i^0$, $\tau_i^a$ to coincide with interval boundaries $t_i$ and $t_{i+1} \, \forall i \in \{0,\dots,N-1\}$. The default values for the number of intervals, state and input polynomial degree and quadrature order are $N=5$, $a=2$, $b=1$ and $Q=3$. 

In order to demonstrate the effectiveness of our method, we showcase the two main features of our proposed method, namely superlinear convergence and the ability to control the accuracy of state and input trajectories, on an optimal control problem with a  chattering solution. While solving optimal control problems is the main focus of our work, the capabilities of the implemented method can be used to solve feasibility problems and complex differential equations as well (such as high index DAEs and differential inclusions). Algorithm \ref{alg:algorithm1} was implemented using the interior point NLP solver Ipopt~\cite{wachter2006implementation} with a relative convergence tolerance set to $10^{-10}$. All tests were performed on a laptop with an Intel\textsuperscript{\textregistered} Core\textsuperscript{\texttrademark} i7-4600U CPU at 2.10\,GHz with 16\,GB of RAM.

\subsection{Fuller problem description}
\label{sec:section4A}
We propose the numerical experiment to be an optimal control problem with a discontinuous solution at non-trivial times in order to underline the capability of our flexible-mesh optimal control solver to capture these discontinuities. Additionally, we will analyse the impact of the desired accuracy on the cost value for the numerically computed solution (which is a lower bound for the exact optimal solution cost).

The chosen problem is a variation of the Fuller problem~\cite{fuller1963study} 
\begin{subequations}
\label{eq:equation12}
\begin{align}
    \min_{p(\cdot), u(\cdot)} \quad & \int_{0}^{T}{p^2(t)}\diff t\label{eq:equation12A}\\
    \textrm{s.t.} \quad & \ddot{p}(t)=u(t), \quad\quad\quad\,\,\,\, \forall t\in[0,T]\label{eq:equation12B}\\
    & u(t)\in[-0.01,0.01], \, \forall t\in[0,T]\label{eq:equation12C}\\
    & p(0)=p(T)=\dot{p}(T)=0, \, \dot{p}(0)=1,\label{eq:equation12D}
\end{align}
\end{subequations}
with $T=300$ seconds, where $p(t)\in\mathbb{R}$ is the position and $u(t)\in\mathbb{R}$ the control input at time $t$. 

\begin{figure}
    \centering
    \includegraphics[width=\columnwidth]{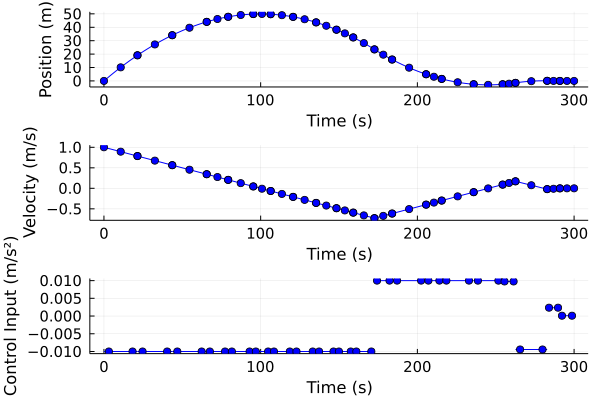}
    \caption{Numerically computed control solution to Fuller problem using flexible meshes with $N=20$ intervals, polynomial degrees $a=2$ , $b=1$, flexibility parameter $\phi=0.5$ and desired accuracy $\epsilon_{tol}=10^{-8}$. Blue dots indicate the location of mesh points.}
    \label{fig:figure1}
\end{figure}

Since the chosen time $T$ was sufficiently large, the optimal control input trajectory has a bang-bang structure with values alternating between $-0.01$ and $0.01$ and then reaching the steady state  with $u(T)=0$. These switching times are difficult to be captured by a numerical solver. Figure~\ref{fig:figure1} displays the state components $p$, $\dot{p}$ and the control input $u$ as obtained by implementing Algorithm~\ref{alg:algorithm1}. Another relevant feature to observe is how the mesh automatically becomes denser in the regions of sudden changes near the switches and coarser where the solution is smoother. 

Another aspect which motivates this choice of  illustrative example is the chattering phenomenon. As  can be observed, instead of getting infinitely many  switches, we only capture a finite number of switches between values that are not all on the input bounds. This behaviour is due to the specified  tolerances, as  explained in Section \ref{sec:section4C} below. The lower the tolerance, the more accurate the numerical solution becomes. 

\subsection{Superlinear convergence for discontinuous solutions}
\label{sec:section4B}
Note that in this problem the solution can be represented by piecewise polynomials of a sufficiently high degree. As a result, when setting $a\geq2$ the minimum in~\eqref{eq:equation9} converges to $1.8729\cdot10^{-15}$. In general, solutions cannot be represented exactly using piecewise polynomials, hence we set the polynomial degree to $a=b=1$ in order to reproduce the convergence behaviour generally encountered for most practical problems. 

\begin{figure}
    \centering
    \includegraphics[width=\columnwidth]{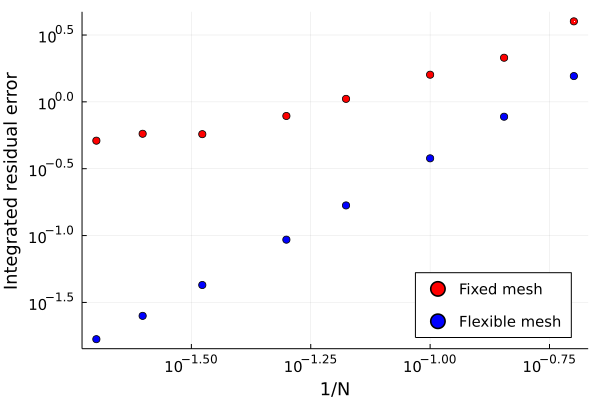}
    \caption{Convergence of the Fuller problem as the mesh is refined with $N$ between 5 and 60 intervals, first order polynomial approximation $a=b=1$ and flexibility parameter $\phi=0.5$.}
    \label{fig:figure2}
\end{figure}

Figure \ref{fig:figure2} presents the performance of our flexible mesh idea and compared to a fixed mesh refinement procedure. The plot shows the variation of the minimum integrated residual attainable for a certain number of mesh intervals $N$. Note the scale is logarithmic and on the horizontal axis is plotted the inverse of $N$. The slope of the blue line is approximately $2$ which is an indicator of superlinear convergence. In contrast with our flexible mesh, the red dots form a line of slope approximately equal to one for low values of $N$ and then start to plateau at a certain value.

\subsection{Between accuracy and optimality}
\label{sec:section4C}
\begin{figure}
    \centering
    \includegraphics[width=\columnwidth]{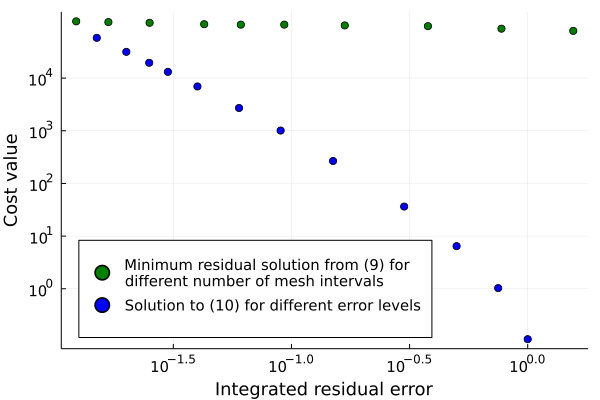}
    \caption{Impact of the desired accuracy on the cost value for the numerically computed solution of the Fuller problem using flexible meshes with $\phi=0.5$. Green dots denote solutions to problem \eqref{eq:equation9} for different $N$ and when the cost is not minimized. Blue dots represent the solutions to problem \eqref{eq:equation11} for fixed parameters $N=60$, $a=b=1$ as the tolerance $\epsilon_{tol}$ is varied.}
    \label{fig:figure3}
\end{figure}

\begin{figure}
    \centering
    \includegraphics[width=\columnwidth]{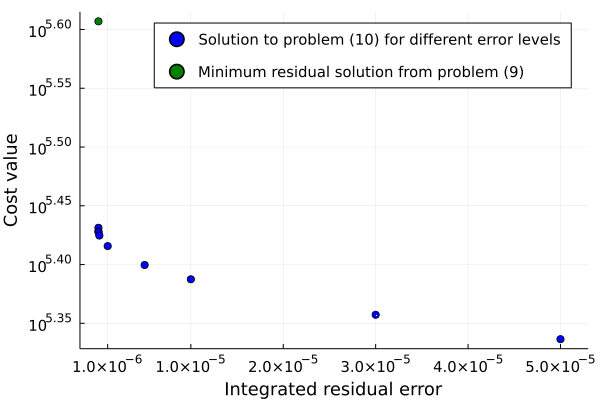}
    \caption{Impact of the desired accuracy on the cost value for the numerically computed solution for the Fuller problem using flexible meshes with $\phi=0.5$. The green dot shows the solution to \eqref{eq:equation9} for $N=20$, $a=2$, $b=1$. Blue dots represent the solutions to problem \eqref{eq:equation11} as the tolerance $\epsilon_{tol}$ is varied.}
    \label{fig:figure4}
\end{figure}

Figure \ref{fig:figure3} presents in green the solution for problem \eqref{eq:equation9} using an increasing number of mesh nodes $N$. Blue denotes the solution for \eqref{eq:equation11} using an increasing tolerance $\epsilon_{tol}$.
The green dots represent what happens to the minimized integrated residual as the number of intervals is changed. As can be observed, the cost does not change significantly as the mesh is refined. However this is not always the case and there can be situations where the solution of problem \eqref{eq:equation9} is very far away from the numerical solution computed from \eqref{eq:equation11}. Such an example is shown in Figure~\ref{fig:figure4}. The blue dots are points on the Pareto front between numerical lower bounds on the optimal cost value and discretization error for a fixed number of intervals $N=20$. 

Having an integrated residual based transcription allows the user to generate an approximate numerical solution using a finite number of discretization intervals $N$ and visualize the impact of different $\epsilon_{tol}$ tolerance values on the computed solution. As one would expect, the objective value increases as the tolerance is decreased, since the solution is captured more accurately by our numerical scheme and a tighter lower bound to the exact solution can be produced.


The Pareto front shown in Figure~\ref{fig:figure4} uses a logarithmic scale for the y-axis and a linear scale for the x-axis. In this example, decreasing the  tolerance $\epsilon_{tol}$ means that more full switches will be captured in the region $t\approx 280$ seconds, thus improving the solution accuracy. However, reducing $\epsilon_{tol}$ will increase the computational time.  As can be noticed, below a certain accuracy threshold the cost value no longer changes significantly, meaning that the gap between the numerical and exact  optimal cost cannot be further reduced without increasing the number of mesh nodes $N$. In general, we aim to find a value for $\epsilon_{tol}$ that can maintain a good balance between solution accuracy and computational time. 

\section{Conclusions and future works}
\label{sec:section5}
Numerically solving optimal control problems for a given discretization mesh involves a trade-off between computational time and solution accuracy. In our example, we can use the change in cost value as a function of the residual tolerance~$\epsilon_{tol}$ as a metric to determine whether a satisfactory solution has been found. While in most transcription methods this trade-off cannot be easily ensured a priori, our proposed algorithm explicitly includes the residual tolerance $\epsilon_{tol}$ as a parameter and is able to construct the Pareto front between the cost value and residual tolerance $\epsilon_{tol}$. This opens up the possibility for early termination in real-time control applications. As discussed, it is important to use a transcription method that relies on error measures that are integrated along the entire solution trajectory to assess convergence, instead of measured error at a finite number of sampled locations. In this way, the error between mesh nodes can be accounted for.   

When discontinuities are present in the solution, numerically approximating state and input trajectories to an acceptable tolerance is especially challenging. Our method proposes the use of a flexible mesh for capturing discontinuities by including time mesh nodes in the decision vector. The efficiency of this method was demonstrated through an illustrative example. 

The implementation of the method is in an early development stage and many improvements are possible in order to demonstrate its full capabilities. Further research could be conducted on improving the mesh refinement process by including early termination procedures, thus increasing the computational efficiency of the overall solution process. Other future work could aim at automatically increasing the quadrature order such that the constraints in~\eqref{eq:equation7} can be removed. Theoretical convergence and performance guarantees also need to be investigated.

\balance
\bibliographystyle{plain}
\bibliography{references.bib}

\addtolength{\textheight}{-12cm}  

\end{document}